\documentstyle[12pt]{article}
\newtheorem{ttt}{Theorem}[section]
\newtheorem{ppp}{Proposition}[section]
\newtheorem{lll}{Lemma}[section]
\newtheorem{ccc}{Corollary}[section]

\begin{document}
\vspace{2cm}
\centerline{\Large\bf ESSENTIALS}\vspace{1cm}
\centerline{\Large\bf OF THE}\vspace{1cm}
\centerline{\Large\bf METHOD OF MAXIMAL SINGULARITIES}\vspace{1cm}

\centerline{\large Aleksandr V. Pukhlikov}\vspace{1cm}
\centerline{Institute for Systems Analysis, Moscow}
\centerline{March 23, 1996}
\section{Introduction}

We start with a very brief description of the principal events
in the history of our subject.

About 1870 M.Noether discovered [N] that
the group of birational automorphisms of the projective plane
$\mathop{\rm Bir}{\bf P}^2$, which is known also as the {\it Cremona
group}
$\mathop{\rm Cr} {\bf P}^2$,
is generated by its subgroup
$\mathop{\rm Aut} {\bf P}^2=\mathop{\rm Aut} {\bf C}^3/{\bf C}^*$ and any
{\it quadratic Cremona transformation} $\tau$,
which in a certain system of homogeneous coordinates can be
written down as
$$\tau\colon(x_0\colon x_1\colon x_2)\to
(x_1x_2\colon x_0x_2\colon x_0x_1).$$
Noether's arguments were as follows. Take any Cremona transformation
$$
\chi\colon{\bf P}^2-\,-\,\to{\bf P}^2.
$$
Then either $\chi$ is a projective isomorphism, or the proper inverse
image of the linear system of lines in ${\bf P}^2$
is a linear system $|\chi|$ of curves of the degree
$n(\chi)\ge2$ with prescribed base points
$a_1,\dots,a_N$, including infinitely near ones.
Let $\nu_1,\dots,\nu_N$ be their multiplicities
with respect to the system $|\chi|$, and assume
that $\nu_1\ge\nu_2\ge\dots\ge\nu_N$. Then, as far as
two lines intersect each other at one point, the
{\it free intersection} of two curves of the system $|\chi|$
(that is, the intersection outside the base locus)
equals 1. So
$$ n^2(\chi)=
\sum\limits^N_{i=1}\nu^2_i+1.$$
Moreover, the curves in $|\chi|$ are rational, and nonsingular outside the base
locus, and so, computing their geometric genus
by their arithmetical one, we get
$$\bigl(n(\chi)-1\bigr)\bigl(n(\chi)-2\bigr)=
\sum\limits^N_{i=1}\nu_i(\nu_i-1).$$
It is easy to derive from these two equalities that $N\ge 3$ and
the first three greatest multiplicities satisfy the
{\it Noether's inequality}
$$\nu_1+\nu_2+\nu_3>n.$$

Now, if $a_1,a_2,a_3$ lie directly on the plane ${\bf P}^2$ (i.e.,
among them there are no infinitely near points), then we can take the
composition
$$\chi\tau\colon{\bf P}^2-\,-\,\to{\bf P}^2,$$
where $\tau$ is a quadratic transformation constructed with respect to
the triplet $(a_1,a_2,a_3)$. Let us prove that $n(\chi\tau)<n(\chi)$.

Indeed, the degree of a curve is the number of the points of intersection
with a line. But the points of intersection of a generic line $L$ and
a curve $C$ from $|\chi\tau|$ are in 1-1 correspondence with the points
of {\it free intersection} of their images $\tau(L)$ and $\tau(C)$. But
$\tau(L)$ is a conic passing through $a_1,a_2,a_3$, whereas $\tau(C)\in|\chi|$.
Thus the intersection number equals
$$2n(\chi)-\nu_1-\nu_2-\nu_3<n.$$

Proceeding in this manner, we ``untwist'' the ``maximal'' triplets until
we come to the case $n(\chi)=1$. Noether's theorem would have been
proved but for the maximal triplets which include infinitely near points,
when Noether's arguments do not work. It took about 30 years to complete
the proof.\newline\newline

The second part of our story begins in the first years of the present
century, when G.Fano made his first attempts to extend two-dimensional
birational methods to three-dimensional varieties [F1,2].

He started with trying to describe birational correspondences of three-dimensional
quartics $V_4\subset {\bf P}^4.$ His choice of the object of study was
really the best: up to this day, the quartic is one of the principle
touchstones for multi-dimensional birational constructions.

Reproducing Noether's arguments, Fano considered a birational correspondence
$\chi\colon V\,-\,-\,\to V'$ between two smooth three-dimensional quartics
and, taking the proper inverse  image $|\chi|$ of the linear system of hyperplane
sections of $V'\subset{\bf P}^4$, came to the following conclusion.
Either $|\chi|$ is cut out on $V$ by hyperplanes,
and then $\chi$ is a projective
isomorphism, or $|\chi|$ is cut out on $V$ by hypersurfaces
of the degree $n(\chi)\ge 2$, and then the base locus $|\chi|$
satisfies certain conditions, similar to the Noether's inequality.
These conditions are now called the {\it Fano's inequalities}. Fano
asserted that, if $n(\chi)\ge 2$, then something like one of
the following two
cases happens:\newline\newline
either there is a point
$x\in V$ such that $\mathop{\rm mult}_x|\chi|>2n$,\newline\newline
or there is a curve $C\subset V$ such that $\mathop{\rm mult}_C|\chi|>n$.
\newline\newline

It should be admitted that Fano never
asserted that {\it only} these two cases are
possible. He gives the following example:
a base point
$x\in V$
and a base line
$L\subset E$,
where $E\subset\tilde V$
is the exceptional divisor of the blowing up
of
$x$,
satisfying the inequality
$$
\mathop{\rm mult}\nolimits_x|\chi|+\mathop{\rm mult}\nolimits_L|\tilde\chi|>3n.
$$
But the general level of understanding and technical
weakness of his time prevented him from giving a
rigorous and complete description of what happens when
$n(\chi)\geq 2$.

Then Fano asserted that none of his conditions can hold.
It is really amazing, that practically all of his arguments
being absolutely invalid, this very assertion is true.
(It is still more amazing that this was very often the case with
Fano: wrong arguments almost always led him to deep and
true conclusions.) For instance, to exclude the
possibility of a curve
$C$ with
$\mathop{\rm mult}\nolimits_C|\chi|>n$,
he employs the arguments of arithmetic genus of
a general surface in
$|\chi|$.
It seems that Fano tried to reproduce Noether's
arguments which used the genus of a curve in
$|\chi|$.
However, Iskovskikh and Manin found out that these
arguments actually lead to no conclusion [IM].

Being sure that  the case
$n(\chi)\geq 2$ is impossible, Fano
formulated one of his most impressive claims:
any birational correspondence between two
non-singular quartics in
${\bf P}^4$ is a projective isomorphism.
In particular, the group  of birational automorphisms
$\mathop{\rm Bir} V=\mathop{\rm Aut} V$
is finite, generically trivial, and
$V$
is nonrational.

In this way Fano did a lot of work in three-dimensional
birational geometry [F3]. He gave a description
(however incomplete and unsubstantiated) of birational correspondences
of three-dimensional cubics, complete intersections
$V_{2\cdot 3}$ in ${\bf P}^5$ and
many other varieties. A lot of his results has not
been completed up to this day.

However, because of the very style of Fano's work, his
numerous mistakes and, generally speaking, incompatibility
of his geometry with the universally adopted standards of
mathematical arguments, his ideas and computations had been
abandoned for a long period of about twenty years.
\newline\newline

In late sixties Yu.I.Manin and V.A.Iskovskikh in Moscow
(after a series of papers on two-dimensional birational
geometry) started their pioneer program in three-dimensional
birational geometry. As a result, in 1970 they developed a
method which was sufficiently strong to prove the Fano's
claim on the three-dimensional quartic [IM]. We shall refer
to this method as to the {\it method of maximal singularities}.
By means of this method Iskovskikh proved later a few more Fano's
claims and corrected some of Fano's mistakes [I]. In seventies-eighties
several students of Iskovskikh -- A.A.Zagorskii, V.G.Sarkisov [S1, S2],
S.L.Tregub [T1, T2], S.I.Khashin [Kh] and the author [P1-P5]
-- had been working in this field, trying to describe birational
correspondences of certain classes of algebraic varieties. Sometimes,
although not very often, their work was successful: the method of
maximal singularities was extended to a number of classes of varieties,
of arbitrary dimension and possibly singular, including a big class
of conic bundles. The well-known Sarkisov program [R1, C1] was born
in the framework of this field, too. At the same time, the method really
works only for varieties of a very small degree. One must admit that
at present we have no good method for studying birational geometry of
multi-dimensional Fano varieties and Fano fibrations.

Nevertheless, the results obtained by means of the method of
maximal singularities can not be proved at present in any other way
(see [K] for an alternative approach in the spirit of $p$-characteristic
tricks).\newline\newline

This paper is an extraction made from the lectures given by the author
during his stay at the University of Warwick in September-December 1995.
Since [IP] has been published, it makes no sense to reproduce all the
details of excluding/untwisting procedures here. At the same time,
[IP] was actually written in 1988. After that the real meaning
of the ``test class'' construction became clearer, and some
new methods of exclusion of maximal singularities appeared [P5,6].
The aim of the present paper is to give an easy introduction to the
method of maximal singularities. We restrict ourselves by
explanation of crucial points only. The principal and most difficult part
of the method -- that is, exclusion of infinitely near maximal
singularities, -- is presented here in the new form, simple and easy.
This version of the method has never been  published before.

I would like to express my gratitude to Professor M.Reid who
invited me to the University of Warwick and arranged my lecture
course
in birational geometry. I am thankful to all the staff of the
Mathematics Institute for hospitality.

I would like to thank Professor V.A.Iskovskikh, who set up before myself
the problem of four-dimensional quintic in 1982 and thus
determined the direction of my work in algebraic geometry.

I am grateful to Professor Yu.I.Manin for constant
and valuable support.

The author was financially supported by International Science
Foundation, grant M9O000, by ISF and Government of Russia,
grant M9O300, and Russian Fundamental Research Fund, grant 93-011-1539.

\section{Maximal singularities of birational maps}

Fix a projective ${\bf Q}$-factorial variety
$V$
with at most terminal singularities
over the field
${\bf C}$
of complex numbers.

\subsection{Test pairs}

{\bf Definition.} A pair
$(W,Y)$, where
$W$
is a projective variety such that
$\mathop{\rm dim} W=\mathop{\rm dim} V$,
$\mathop{\rm codim}\mathop{\rm Sing} W\geq 2$ and
$Y$
is a divisor on
$W$
is said to be a {\it test pair}, if the following conditions hold:
\newline
a)
$|Y|$ is free from fixed components,\newline
b)
there exists a number
$\alpha\in{\bf R}_+$
such that for any
$\beta>\alpha, \beta\in{\bf Q}$
the linear system
$$
|M(Y+\beta K_W)|
$$
is empty for
$M\in{\bf Z}_+, M\beta\in{\bf Z}$
(the adjunction break condition).

The minimal
$\alpha\in{\bf R}_+$
satisfying the condition b)
is said to be the {\it index}
(or {\it threshold})
of the pair
$(W,Y)$.
We shall denote it by
$\alpha(W,Y)$.\newline\newline

Our aim is to study the maps
$$
\chi:V-\,-\,\to W.
$$

\subsubsection*{Examples}

We list the types of varieties, which were among the principal
objects of (more or less succesful) study by
means of the method of maximal singularities during the
last twenty five years:
\newline
\newline

-- a smooth quartic
$V_4\subset{\bf P}^4;$
\newline
\newline

-- a complete intersection
$V_{2\cdot 3}\subset{\bf P}^5$;
\newline
\newline

-- a  singular quartic
$x\in V_4\subset{\bf P}^4$;
\newline
\newline

-- a smooth hypersurface
$V_M\subset{\bf P}^M$;
\newline
\newline

-- a double projective space
$\sigma:V\to{\bf P}^n$
branched over a smooth hypersurface
$Z_{2n}\subset{\bf P}^n$.\newline\newline

Now let us give the principal examples of test pairs,
explaining briefly what do we need them for.
\newline\newline

(${\bf P}^n$, hyperplane) -- to decide whether
$V$ is rational; note that the index of this pair
is equal to $1/(n+1)$;\newline\newline

($\varphi:W\to S$ -- a Fano fibration,
$Y=\varphi^{-1}$(very ample divisor on $S$)) --
to decide whether there are structures of Fano fibrations
on $V$; for instance, take a conic bundle or Del Pezzo
fibration; the index here is obviously zero;\newline\newline

($V,(-MK_V)$) -- to describe the group
$\mathop{\rm Bir} V$ and to give the birational classification
within the same family of Fano varieties.

\subsubsection*{The first step}

Assume that there exists a birational map
$$
\chi:V-\,-\,\to W.
$$
Take the proper inverse image
$|\chi|\subset|D|$
of the linear system
$|Y|$
on
$V$.
Denote by
$\mathop{\rm Bs} |\chi|$
its base subscheme.

This system
$|\chi|$
and this subscheme
$\mathop{\rm Bs} |\chi|$
are the principal objects of our study.

\subsection{The language of discrete valuations}

We shall remind briefly the necessary definitions and facts
about discrete valuations. For more details see [P5,6].

Let
$X$
be a quasi-projective variety.

Denote by
${\cal N}(X)$ the set of {\it geometric}
discrete valuations
$$
\nu:{\bf C}(X)^*\to{\bf Z},
$$
which have a centre on
$X$.
If $X$
is complete, then
${\cal N}(X)$ includes all the geometric discrete valuations.
The centre of a discrete valuation
$\nu\in{\cal N}(X)$
is denoted by
$Z(X,\nu)$.

{\bf Examples.}(1) Let
$D\subset X$ be a prime divisor,
$D\not\subset\mathop{\rm Sing} X$.
Then
$D$
determines a discrete valuation
$$
\nu_D=\mathop{\rm ord}\nolimits_D.
$$
(2) Let
$B\subset X$
be an irreducible cycle,
$B\not\subset\mathop{\rm Sing} X$.
Then
$B$
determines a discrete valuation:
$$
\nu_B(f)=
\mathop{\rm mult}\nolimits_B(f)_0-
\mathop{\rm mult}\nolimits_B(f)_{\infty}.
$$

Note that if
$\sigma_B:X(B)\to X$
is the blowing up of
$B$,
$E(B)=\sigma^{-1}_B(B)$
is the exceptional divisor, then
$$
\nu_B=\nu_{E(B)},
$$
where
${\bf C}(X)$
and
${\bf C}(X(B))$
are naturally identified.

{\bf Definition.} Let
$\nu\in{\cal N}(X)$
be a discrete valuation. A triplet
$(\tilde X,\varphi,H)$,
where
$\varphi:\tilde X\to X$
is a birational morphism,
$H\not\subset\mathop{\rm Sing}\tilde X$
is a prime divisor, is called a {\it realization} of
$\nu$,
if
$\nu=\nu_H$.

\subsubsection*{Multiplicities}

Let
$\nu\in{\cal N}(X),{\cal J}\subset{\cal O}_X$
be a sheaf of ideals.

{\bf Definition.} The multiplicity of
${\cal J}$ at
$\nu$
equals
$$
\nu({\cal J})=
\mathop{\rm mult}\nolimits_H\varphi^*{\cal J},
$$
where
$(\tilde X,\varphi,H)$
is a realization of
$\nu$.

Let
$|\lambda|\subset|D|$
be a linear system of Cartier divisors,
${\cal L}(|\lambda|)\subset{\cal O}_X(D)$
be the subsheaf generated by the global sections in
$|\lambda|$. Set
$$
{\cal J}(|\lambda|)={\cal L}(|\lambda|)\otimes{\cal O}_X(-D)\subset{\cal O}_X.
$$
Obviously,
${\cal J}(|\lambda|)$
is the ideal sheaf of the base subscheme
$\mathop{\rm Bs}|\lambda|$.

{\bf Definition.} The multiplicity of
$|\lambda|$ at
$\nu$
equals
$$
\nu(|\lambda|)=
\nu({\cal J}(|\lambda|)).
$$

Now let
$X$
be
(${\bf Q}$-)Gorenstein,
$\pi:X_1\to X$
be a resolution. Then
$$
K_{X_1}=\pi^*K_X+\sum_i d_iE_i
$$
for some prime divisors
$E_i\subset X_1$.
Take a realization
$(\tilde X,\varphi,H)$,
of
$\nu\in{\cal N}(X_1)={\cal N}(X)$.
Then we get an inclusion
$$
\varphi^*\omega_{X_1}\hookrightarrow\omega_{\tilde X},
$$
and the following ideal sheaf on
$\tilde X$:
$$
K(X_1,\varphi)=
\omega^*_{\tilde X}\otimes
\varphi^*\omega_{X_1}
\hookrightarrow{\cal O}_{\tilde X}.
$$

{\bf Definition.} The {\it canonical multiplicity (discrepancy)}
of
$\nu$
is equal to
$d_i$,
if
$\nu=\nu_{E_i}$,
and to
$$
K(X,\nu)=
\mathop{\rm mult}\nolimits_HK(X_1,\varphi)+
\sum_id_i\nu(E_i),
$$
otherwise.

{\bf Example.} Let
$B\subset X$,
$B\not\subset\mathop{\rm Sing} X$
be an irreducible cycle of codimension$\geq 2$. Then
$$
\nu_B({\cal J})=\mathop{\rm mult}\nolimits_B{\cal J},
$$
$$
\nu_B(|\lambda|)=\mathop{\rm mult}\nolimits_B|\lambda|,
$$
$$
K(X,\nu_B)=\mathop{\rm codim} B-1.
$$

\subsection{Maximal singularities}

Let us return to our variety
$V$ and birational map
$\chi:V-\,-\,\to W$.
Denote by
$n(\chi)$
the index
(threshold) of the pair
$(V,D)$.

{\bf Definition.}
A discrete valuation
$\nu\in{\cal N}(V)$
is said to be a {\it maximal singularity} of
$\chi$,
if the following inequality holds:
$$
\nu(|\chi|)>n(\chi)K(V,\nu).
$$

\begin{ttt}
\label{1}
Either
$\alpha(V,D)\leq
\alpha(W,Y),$
or
$\chi$
has a maximal singularity.
\end{ttt}

{\bf Proof:} see [P5,6]. It is actually so easy that
can be left to the reader as an exercise. The idea of the proof
can be found in any paper concerned with these problems
(for instance, [IM,I,P1,IP]). Just keep in mind that
the proof should not depend upon resolution of singularities.

{\bf Example.}
Let
$V$
be smooth with
$\mathop{\rm Pic}\cong{\bf Z} K_V$
and assume that the anticanonical system
$|-K_V|$ is free. Then
$$
|\chi|\subset|-n(\chi)K_V|,
$$
and for a birational automorphism
$\chi\in\mathop{\rm Bir} V$
either
$n(\chi)=1$,
or
$\chi$
has a maximal singularity.

\subsubsection*{Maximal cycles}

Let
$V$
be non-singular.

{\bf Definition.}
An irreducible cycle
$B\not\subset\mathop{\rm Sing} V$
of codimension$\geq 2$
is said to be a {\it maximal cycle},
if
$\nu_B$
is a maximal singularity. Explicitly:
$$
\mathop{\rm mult}\nolimits_B|\chi|>
n(\chi)(\mathop{\rm codim} B-1).
$$

{\bf Definition.}
A maximal singularity
$\nu\in{\cal N}(V)$
is said to be
{\it infinitely near},
if it is not a maximal cycle.

{\bf Remark.} The meaning of these two definitions
is to separate ``shallow'' maximal singularities,
which are not very far from the ``ground''
$V$,
and ``deep'' ones, which take a lot of blow-ups to get to
them. When
$V$ is singular, these definitions should be
modified slightly by adding some valuations sitting at
singularities (see [P3,6]).

\section{The untwisting scheme}

Assume that
$\alpha(v,D)>\alpha(W,Y)$.
Then
$\chi$
has a maximal singularity.
The untwisting scheme gives an idea of simplifying
$\chi$ according to its maximal singularities.

\subsection{Basic Conjecture}

We say that
$V$
satisfies the Basic conjecture, if for any
$\chi:V-\,-\,\to W$
in the hypothesis of Theorem \ref{1} we can replace
the words ``maximal singularity'' by the words
``maximal cycle'': if
$$
\alpha(v,D)>\alpha(W,Y),
$$
then
$\chi$
has a maximal cycle.

\subsection{Excluding maximal cycles}

Assume that
$V$
satisfies the Basic conjecture. Then the first thing to be done
is to describe all the cycles
$B\subset V$
which can occur as maximal. In other words, all the cycles
$B$
such that for some
$D\in\mathop{\rm Pic} V$
$$
|D-\nu B|
$$
is free from fixed components for some
$$
\nu>(\mathop{\rm codim} B-1)\alpha(V,D).
$$

\subsection{Untwisting maps}

The second step of the scheme is to construct an automorphism
$\tau_B\in\mathop{\rm Bir} V$
for each
$B$
singled out at the previous step. The cycle
$B$
should be maximal for
$\tau_B$.

\subsection{Untwisting}

If
$B$ is maximal for
$\chi:V-\,-\,\to W$,
take the composition
$$
\chi\circ\tau_B:V-\,-\,\to W.
$$
It {\it should}
turn out that
$$
n(\chi\circ\tau_B)<n(\chi).
$$
Iterating, we come to a sequence of cycles
$B_1,\dots,B_k$
such that
$$
n(\chi\circ\tau_{B_1}\circ\dots\circ\tau_{B_k})\leq
\alpha(W,Y).
$$

\subsection{Birationally rigid varieties}

Informally speaking,
$V$
is birationally rigid, if the untwisting scheme works on it.

{\bf Definition.}
$V$
is said to be
{\it birationally rigid}, if for any test pair
$(W,Y)$
and any map
$\chi:V-\,-\,\to W$
there exists
$\chi^*\in\mathop{\rm Bir} V$
such that
$$
n(\chi\circ\chi^*)\leq
\alpha(W,Y).
$$
If, moreover,
$\mathop{\rm Bir} V=\mathop{\rm Aut} V$,
then
$V$
is said to be
{\it birationally superrigid}.

{\bf Remark.} The untwisting scheme, when it works, gives not
only the fact of birational rigidity, but also a set of natural
generators of the group
$\mathop{\rm Bir} V$ --
that is, the maps
$\tau_B$.

\begin{ppp}
\label{p1}
Assume that
$V$ is birationally rigid and
$\mathop{\rm Pic} V\cong{\bf Z}$.
Then
$V$
has no structures of Fano fibrations.
\end{ppp}

{\bf Proof.} Assume that there is a map
$$
\chi:V-\,-\,\to W,
$$
where
$p:W\to S$
is a Fano fibration. Take
$Y$
to be
$p^{-1}(Q)$,
where
$Q\subset S$
is
a very ample divisor. Then
$$
n(\chi\circ\chi^*)=0
$$
for some
$\chi^*\in\mathop{\rm Bir} V$,
so that
$$
|\chi\circ\chi^*|\subset
|-n(\chi\circ\chi^*)K_V|=|0|.
$$
Contradiction.

In particular,
$V$ has no structures of a conic bundle or a Del Pezzo
fibration.

\section{Excluding maximal cycles}

We show by example how it is to be done. A lot of other
examples can be found in the original papers [IM, I, P1-6, IP].

\subsection{Double spaces}

Let
$\sigma:V\to{\bf P}^m\supset W_{2m}$
be a smooth double space of the index 1,
$m\geq 3$,
branched over a smooth hypersurface
$W$ of the degree
$2m$.
Let
$|\chi|\subset|-nK_V|$
be a system free from fixed components.

\begin{ttt}
\label{t2}
$|\chi|$
has no maximal cycles.
\end{ttt}

\begin{ccc}
Modulo Basic conjecture
$V$ is superrigid.
\end{ccc}

{\bf Proof.} It breaks into two parts: we exclude maximal
points and maximal cycles of positive dimension separately.

\subsection{Points}

Obviously, a point
$x\in V$
cannot be maximal: take a plane
$\bar P\ni\bar x=\sigma(x),$
then
$P=\sigma^{-1}(\bar P)$
is a nonsingular surface,
$|\chi|_P$
has no fixed curves, so that for any
$D_1,D_2\in|\chi|_P$
$$
(D_1\cdot D_2)=2n^2.
$$
But
$\mathop{\rm mult}\nolimits_xD_i>2n:$ a contradiction.

\subsection{Curves}

\begin{ppp}
For any curve
$C\subset V$
$$
\mathop{\rm mult}\nolimits_C|\chi|\leq n.
$$
\end{ppp}
Obviously, our theorem is an immediate consequence of this fact.

{\bf Proof of the Proposition.}
Let us consider the following three cases:\newline\newline
(1) $C=\sigma^{-1}(\bar C), \bar C\not\subset W$; \newline\newline
(2) $\bar C=\sigma(C)\subset W$; \newline\newline
(3) $\sigma:C\to\bar C$ is birational, $\bar C\not\subset W$.

\subsubsection*{The easy first case.}

Take a generic line
$\bar L$
intersecting
$\bar C$,
$L=\sigma^{-1}(\bar L)$
is a smooth curve. The linear series
$$
|\chi|\Bigr|_L
$$
is of the degree
$2n$
and has $\geq 2$
points$\in\sigma^{-1}(\bar L\cap\bar C)$
of the multiplicity
$\mathop{\rm mult}\nolimits_C|\chi|$.

\subsubsection*{The second case, not very difficult.}

Take a generic point
$x\in{\bf P}^m$
and the cone
$Z(x)$
over
$\bar C$
with the vertex
$x$.
Then
$Z(x)\cap W=\bar C\cup\bar R(x)$,
where the residual curve
$\bar R(x)$
intersects
$\bar C$
at
$\mathop{\rm deg} R(x)$
different points
(see [P5]). Let
$R(x)$
be the curve
$\sigma^{-1}(\bar R(x))$,
then
$\sigma:R(x)\to\bar R(x)$
is an isomorphism, and
$$
|\chi|\Bigr|_{R(x)}
$$
is a linear series of the degree
$n\mathop{\rm deg}\bar R(x)$
which has
$\mathop{\rm deg} \bar R(x)$
base points of the multiplicity
$\mathop{\rm mult}\nolimits_C|\chi|$.

\subsubsection*{The non-trivial third case.}

Again take a generic point
$x\in{\bf P}^m$
and the cone
$Z(x)$
over
$\bar C$
with the vertex
$x$.
Let
$$
\varphi:X\to{\bf P}^m
$$
be the blowing up of
$x$
with the exceptional divisor
$E$, so that the projection
$$
\pi:X\to{\bf P}^{m-1}=E
$$
is a regular map,
$X$
being a
${\bf P}^1$-bundle
over
$E$.
Let
$$
\alpha:Q\to\bar C
$$
be the desingularization of
$\bar C$,
$$
\bar S=Q\mathop{\times}\nolimits_{\pi(\bar C)}X
$$
be a
${\bf P}^1$-bundle over $Q$.
Obviously,
$\mathop{\rm Num}\bar S=A^1(S)={\bf Z} f\oplus{\bf Z} e$,
where
$f$
is the class of a fiber and
$e$
is the class of the exceptional section coming from the
vertex of the cone. Obviously,
$f^2=0,(f\cdot e)=1,e^2=-d$,
where
$d=\mathop{\rm deg}\bar C=\mathop{\rm deg}\pi(\bar C).$
Let
$h$
be the class of a hyperplane section,
$h^2=d$,
so that
$h=e+df$.

Denote by
$\tilde C$
the inverse image of
$\bar C$ on
$\bar S$.
Obviously, its class
$\tilde c$
equals
$h$.

For a generic point
$x$
the set
$\sigma^{-1}(Z(x))\cap\mathop{\rm Bs}|\chi|$
contains at most two curves:
$C$
itself and the other component of
$\sigma^{-1}(\bar C)$;
moreover, the inverse image
$\bar W$
of $W$ on
$\bar S$
is a non-singular curve.

Now let us take the surface
$S=\bar S\times_{Z(x)}V$,
that is, the double cover of
$\bar S$
with the smooth branch divisor
$\bar W$.
Denote the image of $C$
on
$S$
by
$C$ again, the other component of
$\sigma^{-1}(\bar C)$
on
$S$
by
$C^*$.
The inverse image of the linear system
$|\chi|$
on
$S$
has at most two fixed components
$C,C^*$
of the multiplicities
$\nu,\nu^*$,
respectively. Therefore the system
$|nh-\nu c-\nu c^*|$
is free from fixed components, and we get the following inequalities:
$$
\Bigl(
(nh-\nu c-\nu^*c^*)\cdot c
\Bigr)\geq 0,
$$
$$
\Bigl(
(nh-\nu c-\nu^*c^*)\cdot c^*
\Bigr)\geq 0.
$$
It is easy to compute the multiplication table for the classes
$h,c$ and $c^*$.
The only necessary intersection number
(the other are obvious) is
$$
(c\cdot c^*)_S=
\frac12(\tilde c\cdot\bar w)_{\bar S}=md.
$$
Now we get the following system of linear inequalities:
$$
(n-\nu^*)+(m-1)(\nu-\nu^*)\geq 0,
$$
$$
(n-\nu)+(m-1)(\nu^*-\nu)\geq 0.
$$
If, for instance,
$\nu\geq\nu^*$,
then by the second inequality
$\nu\leq n$.
By symmetry we are done.
Q.E.D.

\subsection{The general idea of exclusion}

It is very simple: to construct a sufficiently big family of
curves or surfaces intersecting the cycle being excluded at as
many points as possible (or containing it) and, at the same time,
having as small ``degree'' as possible.

Then we restrict our linear system to such a curve or surface
and get a contradiction.

\subsection{What do we know about maximal cycles}

They do not exist:\newline\newline
for smooth hypersurfaces of the degree $M$ in
${\bf P}^M, M\geq 4$[P5]; \newline\newline
for smooth double spaces
$V_2\to{\bf P}^M\supset W_{2M}, M\geq 3$:
[I] for $M=3$, [P2] for $M\geq 4$, see also [IP]
(and for slightly singular as well [P6]);\newline\newline
for smooth double quadrics
$V_4\to Q_2\subset{\bf P}^{M+1}$,
branched over
$Q_2\cap W_{2M-2}, M\geq 4$ [P2], see also [IP].\newline\newline

For a singular quartic
$V_4\subset{\bf P}^4$
with a unique double singular point
$x$ there can be only 25
maximal cycles:
either
$x$ itself, or one of 24 lines on
$V$, containing $x$. Moreover, a
maximal cycle is always unique
[P3].\newline\newline

For a double quadric
$\sigma:V_4\to Q_2\subset{\bf P}^4$,
branched over
$Q_2\cap W_4$,
there can be at most one maximal cycle,
that is, a line
$L\subset V,(L\cdot K_V)=-1,\sigma(L)\not\subset W_4$
[I,IP].\newline\newline

For a complete intersection
$V=V_{2\cdot 3}=Q_2\cap Q_3\subset{\bf P}^5$
a maximal cycle
$B$ is a curve:
either a line
$L$,
or a smooth conic
$Y$
such that the unique plane
$P(Y)\supset Y$
is contained in the quadric
$Q_2$.
Moreover, there can be at most two maximal curves, and if
there are exactly two maximal curves, then they are lines
$L_1$ and $L_2$ such that the unique plane
$P(L_1\cup L_2)\supset(L_1\cup L_2)$
is contained in
$Q_2$ [I,IP].

\section{Untwisting maximal cycles}

We give a simple example of untwisting
(probably the simplest one): the untwisting procedure
for the maximal singular point
$x\in V_4\subset{\bf P}^4$
on a singular quartic $V$ [P3].

\subsection{Construction of the untwisting map}

Let
$\pi:V\backslash\{x\}\to{\bf P}^3$
be the projection from
$x$,
$\mathop{\rm deg}\pi=2$.
Then the untwisting map
$\tau:V-\,-\,\to V$
permutes the points in the fibers of
$\pi$.

Let
$\sigma:V_0\to V$
be the blowing up of
$x$,
$E=\sigma^{-1}(x)\cong{\bf P}^1\times{\bf P}^1$
be the exceptional divisor,
$L_i,i=1,\dots,24$,
be the proper inverse images of lines on
$V$,
passing through
$x$.

\begin{lll}
$\tau$ extends to an automorphism
of
$$
V_0\setminus\mathop{\bigcup}\limits_{1\leq i\leq 24}L_i.
$$
Its action on
$\mathop{\rm Pic} V_0={\bf Z} h\oplus{\bf Z} e$
is given by the following relations:
$$
\tau^* h=3h-4e,
$$
$$
\tau^* e=2h-3e.
$$
\end{lll}

{\bf Proof.}
$\pi$ extends to a morphism
$V\to{\bf P}^3$
of the degree 2. It is not well-defined only on
the one-dimensional fibers, which are exactly the 24 lines
$L_i$.

Thus
$\tau$
is an automorphism of the complement of a set of codimension 2,
so that its action
$\tau^*$
on
$\mathop{\rm Pic} V$
is well-defined.

Obviously, for any plane
$P\subset{\bf P}^3$
its inverse image
$\pi^{-1}(P)$
represents an invariant class,
$$
\tau^*(h-e)=h-e.
$$
Furthermore,
$\pi(E)$
is a quadric in
${\bf P}^3$,
$\pi(H)$
is a quartic in
${\bf P}^3$,
where
$H\subset V$
is a hyperplane section disjoint from
$E$. Thus
$$
e+\tau^* e=2(h-e),
$$
$$
h+\tau^* h=4(h-e).
$$
Q.E.D.

\subsection{Untwisting}

Let
$\chi:V-\,-\,\to W$
be our birational map. We define the number
$\nu_x(\chi)\in{\bf Z}_+$
in the following way: the class of the proper inverse image of the
linear system
$|\chi|$
on
$V_0$
is
$$
n(\chi)h-\nu_x(\chi)e.
$$
The condition ``the singular point $x$
is maximal for
$|\chi|$'' means  that
$$
\nu_x(\chi)>n(\chi).
$$
Now consider the composition
$\chi\circ\tau:V-\,-\,\to W$.

\begin{lll}
(i)
$n(\chi\circ\tau)=3n(\chi)-2\nu_x(\chi).$\newline
(ii) $\nu_x(\chi\circ\tau)=4n(\chi)-3\nu_x(\chi).$
\end{lll}

{\bf Proof.}
Since
$\tau$
is an automorphism in codimension 1, we can write down
$$
n(\chi\circ\tau)h-\nu_x(\chi\circ\tau)e=
\tau^*\Bigl(
n(\chi)h-\nu_x(\chi)e\Bigr).
$$
Applying the formulae, obtained in the previous section,
we get our Lemma.

Now if
$x$
is a maximal point for
$\chi$, then
$\nu_x(\chi)>n(\chi)$,
so that
$n(\chi\circ\tau)<n(\chi)$
and
$\nu_x(\chi\circ\tau)<n(\chi\circ\tau)$,
and
$x$ is no longer a maximal cycle.

The maximal cycle
$x$
is untwisted.

\section{Infinitely near maximal singularities. I}

The techniques necessary to exclude infinitely near
maximal singularities is developed.

\subsection{Resolution}

Let
$X$
be any quasi-projective variety,
$\nu\in{\cal N}(X)$
be a discrete valuation,
$B=Z(X,\nu)\not\subset\mathop{\rm Sing} X, \mathop{\rm codim} B\geq 2$.

\begin{ppp}
Either
$\nu=\nu_B$,
or for the blow up
$$
\sigma_B:X(B)\to X,
$$
$$
E(B)=\sigma^{-1}_B(B),
$$
we get:
$\nu\in{\cal N}(X(B)), Z(X(B),\nu)\subset E(B)$
is an irreducible cycle of codimension$\geq 2$ and
$$
\sigma_B(Z(X(B),\nu))=B.
$$
\end{ppp}

{\bf Proof:} easy.

Consider the sequence of blow ups
$$
\varphi_{i,i-1}:X_i\to X_{i-1},
$$
$i\geq 1$, where
$X_0=X, \varphi_{i,i-1}$
blows up the cycle
$B_{i-1}=Z(X_{i-1},\nu)$
of codimension$\geq 2$,
$E_i=\varphi^{-1}_{i,i-1}(B_{i-1})\subset X_i$.

Set also for
$i>j$
$$
\varphi_{i,j}=\varphi_{j+1,j}\circ\dots\circ
\varphi_{i,i-1}:X_i\to X_j,
$$
$$
\varphi_{i,i}=\mathop{\rm id}\nolimits_{X_i}.
$$
For any cycle
$(\dots)$
we denote its proper inverse image on
$X_i$
by adding the upper index $i$:
$(\dots)^i$.

Note that
$\varphi_{i,j}(B_i)=B_j$
for $i\geq j$.

Note also that although all the $X$'s
are possibly singular,
$B_i\not\subset\mathop{\rm Sing} X_i$
for all $i$.

\begin{ppp}
This sequence is finite: for some
$K\in{\bf Z}_+$
the triplet
$(X_K,\varphi_{K,0},E_K)$
realizes
$\nu$,
$\nu=\nu_{E_K}$.
\end{ppp}

{\bf Proof:} see [P6] or prove it yourself  (it is easy).

{\bf Definition.}
The sequence
$\varphi_{i,i-1}, i=1,\dots,K$,
is said to be the
{\it resolution} of
the discrete valuation
$\nu$
(with respect to the model $X$).

\subsection{The graph structure}

{\bf Definition.}
For
$\mu,\nu\in {\cal N}(X)$
set
$$
\mu\mathop{\leq}\limits_X\nu
$$
if the resolution of
$\mu$ is
a part of the resolution of
$\nu$.

In other words, for some
$L\leq K$
$$
(X_L,\varphi_{L,0},E_L)
$$
is a realization of
$\mu$.

{\bf Definition.}
We define an oriented graph structure on
${\cal N}(X)$, drawing an arrow from
$\nu$
to
$\mu$,
$$
\nu\mathop{\to}\limits_X\mu,
$$
if
$\mu\mathop{\leq}\limits_X \nu$
and
$B_{K-1}\subset E^{K-1}_L$.

Denote by
$P(\nu,\mu)$
the set of all paths from
$\nu$
to
$\mu$
in
${\cal N}(X)$,
which is non-empty if and only if
$\nu\mathop{\geq}\limits_X\mu.$
Set
$$
p(\nu,\mu)=|P(\nu,\mu)|,
$$
if
$\nu\neq\mu$,
and
$p(\nu,\nu)=1$.
Set
${\cal N}(X,\nu)$
to be the subgraph of
${\cal N}(X)$
with the set of vertices smaller (or equal) than $\nu$.

\subsection{Intersections, degrees and multiplicities}

Let
$B\subset X, B\not\subset\mathop{\rm Sing} X$
be an irreducible cycle of codimension$\geq 2$,
$\sigma_B:X(B)\to X$
be, as usual, its blowing up,
$E(B)=\sigma^{-1}_B(B)$
be the exceptional divisor. Let
$$
Z=\sum m_iZ_i,
$$
$$
Z_i\subset E(B)
$$
be a $k$-cycle,
$k\geq\mathop{\rm dim} B$.
We define the
{\it degree}
of $Z$ as
$$
\mathop{\rm deg} Z=
\sum_im_i\mathop{\rm deg}\left(
Z_i\bigcap\sigma^{-1}_B(b)
\right),
$$
where
$b\in B$
is a generic point,
$\sigma^{-1}_B(b)\cong{\bf P}^{\mathop{\rm codim} B-1}$
and the right-hand side degree is the ordinary degree in
the projective space.

Note that
$\mathop{\rm deg} Z_i=0$
if and only if
$\sigma_B(Z_i)$
is a proper closed subset of
$B$.

Our computations will be based upon the following statement.

Let
$D$
and
$Q$
be two different prime Weyl divisors on
$X$, $D^B$ and $Q^B$
be their proper inverse images on
$X(B)$.

\begin{lll}
(i)
Assume that
$\mathop{\rm codim} B\geq 3$.
Then
$$
D^B\bullet Q^B=
(D\bullet Q)^B+Z,
$$
where
$\bullet$
stands for the cycle of the scheme-theoretic intersection,
$\mathop{\rm Supp} Z\subset E(B)$,
and
$$
\mathop{\rm mult}\nolimits_B(D\bullet Q)=
(\mathop{\rm mult}\nolimits_BD)(\mathop{\rm mult}\nolimits_BQ)+\mathop{\rm deg} Z.
$$
(ii) Assume that
$\mathop{\rm codim} B=2$. Then
$$
D^B\bullet Q^B=Z+Z_1,
$$
where
$\mathop{\rm Supp} Z\subset E(B),
\mathop{\rm Supp}\sigma_B(Z_1)$
does not contain $B$,
and
$$
D\bullet Q=
\left[
(\mathop{\rm mult}\nolimits_BD)(\mathop{\rm mult}\nolimits_BQ)+\mathop{\rm deg} Z
\right]
B+
(\sigma_B)_*Z_1.
$$
\end{lll}

{\bf Proof.} Let
$b\in B$
be a generic point,
$S\ni b$
be a germ of a non-singular surface in general position
with $B$,
$S^B$
its proper inverse image on
$X(B)$.
We get an elementary two-dimensional problem:
to compute the intersection number of two different
irreducible curves at a smooth point on a surface
in terms of its blowing up. This is easy.
Q.E.D.

\subsubsection*{Multiplicities in terms of the resolution}

We divide the resolution
$\varphi_{i,i-1}:X_i\to X_{i-1}$
into
the {\it lower part},
$i=1,\dots,L\leq K$,
when
$\mathop{\rm codim} B_{i-1}\geq 3$,
and the {\it upper part},
$i=L+1,\dots,K$,
when
$\mathop{\rm codim} B_{i-1}=2$.
It may occur that
$L=K$
and the upper part is empty.

Let
$|\lambda|$
be a linear system on $X$ with no fixed components,
$|\lambda|^j$
its proper inverse image on
$X_j$.
Set
$$
\nu_j=\mathop{\rm mult}\nolimits_{B_{j-1}}|\lambda|^{j-1}.
$$
Obviously,
$$
\nu_{E_j}(|\lambda|)=
\sum^j_{i=1}p(\nu_{E_j},\nu_{E_i})\nu_i
$$
and
$$
K(X,\nu_{E_j})=
\sum^j_{i=1}p(\nu_{E_j},\nu_{E_i})(\mathop{\rm codim} B_{i-1}-1).
$$
For simplicity of notations we write
$$
i\to j
$$
instead of
$$
\nu_{E_i}\mathop{\to}\limits_X\nu_{E_j}.
$$

Now everything is ready for the principal step of the theory.

\section{Infinitely near maximal singularities. II. \protect\\
The principal computation.}

We prove the crucial inequalities which enable us to exclude
infinitely near maximal singularities for the cases of low degree.

\subsection{Counting multiplicities}

Let
$D_1,D_2\in|\lambda|$
be two different generic divisors. We define a sequence of
codimension 2 cycles on $X_i$'s setting
$$
\begin{array}{l}
\displaystyle
D_1\bullet D_2=Z_0,\\
\displaystyle
D^1_1\bullet D^2_2=Z^1_0+Z_1,\\
\displaystyle
\dots,\\
\displaystyle
D^i_1\bullet D^i_2=
(D^{i-1}_1\bullet D^{i-1}_2)^i+Z_i,\\
\displaystyle
\dots,
\end{array}
$$
where
$Z_i\subset E_i$.
Thus for any
$i\leq L$
we get
$$
D^i_1\bullet D^i_2=
Z^i_0+Z^i_1+\dots+Z^i_{i-1}+Z_i.
$$
For any
$j>i,j\leq L$
set
$$
m_{i,j}=\mathop{\rm mult}\nolimits_{B_{j-1}}(Z^{j-1}_i)
$$
(the multiplicity of an irreducible cycle along a smaller
cycle is understood in the usual sense; for an arbitrary cycle
we extend the multiplicity by linearity).

\subsubsection*{The crucial point}

\begin{lll}
If
$m_{i,j}>0$,
then
$i\to j$.
\end{lll}

{\bf Proof.} If
$m_{i,j}>0$,
then some component of
$Z^{j-1}_i$
contains
$B_{j-1}$.
But
$Z^{j-1}_i\subset E^{j-1}_i$.
Q.E.D.

\subsubsection*{Degree and multiplicity}

Set
$d_i=\mathop{\rm deg} Z_i$.

\begin{lll}
For any
$i\geq 1, j\leq L$
we have
$$
m_{i,j}\leq d_i.
$$
\end{lll}

{\bf Proof.}
The cycles
$B_a$
are non-singular at their generic points. But since
$\varphi_{a,b}:B_a\to B_b$
is surjective, we can count multiplicities at generic points.
Now the multiplicities are non-increasing with respect to
blowing up of a non-singular cycle, so we are reduced to the
obvious case of a hypersurface in a projective space.
Q.E.D.

\subsubsection*{The very computation}

We get the following system of equalities:
$$
\begin{array}{l}
\displaystyle
\nu^2_1+d_1=m_{0,1},\\
\displaystyle
\nu^2_2+d_2=m_{0,2}+m_{1,2},\\
\displaystyle
\vdots\\
\displaystyle
\nu^2_i+d_i=m_{0,i}+\dots+m_{i,i-1},\\
\displaystyle
\vdots\\
\displaystyle
\nu^2_L+d_L=m_{0,L}+\dots+m_{L-1,L}.
\end{array}
$$
Now
$$
d_L\geq
\sum^K_{i=L+1}
\nu^2_i\mathop{\rm deg}(\varphi_{i-1,L})_*B_{i-1}\geq
\sum^K_{i=L+1}\nu^2_i.
$$

{\bf Definition.}
A function
$a:\{1,\dots,L\}\to{\bf R}_+$
is said to be
{\it compatible with the graph structure}, if
$$
a(i)\geq\sum_{j\to i}a(j)
$$
for any
$i=1,\dots,L$.

{\bf Examples:}
$a(i)=p(L,i)$,
$a(i)=p(K,i)$.

\begin{ttt}
Let
$a(\cdot)$
be any compatible function. Then
$$
\sum^L_{i=1}a(i)m_{0,i}\geq
\sum^L_{i=1}a(i)\nu^2_i+
a(L)\sum^K_{i=L+1}\nu^2_i.
$$
\end{ttt}

{\bf Proof.}
Multiply the $i$-th equality by
$a(i)$
and put them all together: in the right-hand part for any
$i\geq 1$
we get the expression
$$
\sum_{j\geq i+1}a(j)m_{i,j}=
\sum_{
\begin{array}{c}
\scriptstyle
j\geq i+1\\
\scriptstyle
m_{i,j}\neq 0
\end{array}
}a(j)m_{i,j}\leq
d_i\sum_{j\to i}a(j)\leq
a(i)d_i.
$$
In the left-hand part for any
$i\geq 1$
we get
$$
a(i)d_i.
$$
So we can throw away all the
$m_{i,*},i\geq 1$,
from the right-hand part, and all the
$d_i, i\geq 1$,
from the left-hand part, replacing
$=$ by $\leq$.
Q.E.D.

\begin{ccc}
Set
$m=m_{0,1}=\mathop{\rm mult}\nolimits_{B_0}(D_1\bullet D_2), D_i\in |\chi|$.
Then
$$
m\left(
\sum^L_{i=1}a(i)
\right)\geq
\sum^L_{i=1}a(i)\nu^2_i+
a(L)\sum^K_{i=L+1}\nu^2_i.
$$
\end{ccc}

\subsection{Applications}

\begin{ccc}
Set
$r_i=p(K,i)$. Then
$$
m\left(
\sum^L_{i=1}r_i
\right)
\geq
\sum^K_{i=1}r_i\nu^2_i.
$$
\end{ccc}

{\bf Proof:}
for
$i\geq L+1$ obviously
$r_i\leq r_L$.
Q.E.D.

\begin{ccc}
{\rm (Iskovskikh and Manin [IM]).}
Let
$\mathop{\rm dim} V=3,\nu\in{\cal N}(V)$
be a maximal singularity such that
$Z(V,\nu)=x$ - a smooth point,
$m=\mathop{\rm mult}\nolimits_xC,$
where the curve $C=(D_1\bullet D_2)$ is the intersection of
two generic divisors from $|\chi|$,
$n=n(\chi)$
and  assume
$|-K_V|$ to be free. Then
$$
m\left(
\sum^L_{i=1}r_i
\right)
\left(
\sum^K_{i=1}r_i
\right)>
n^2\left(
2
\sum^L_{i=1}r_i
+
\sum^K_{i=L+1}r_i
\right)^2.
$$
In particular,
$m>4n^2.$
\end{ccc}

{\bf Proof.}
It follows immediately from the fact that
$\nu$
is a maximal singularity and
the previous Corollary.
Let us prove the last statement.
Denoting
$$
\sum^L_{i=1}r_i,
\sum^K_{i=L+1}r_i
$$
by
$\Sigma_0,\Sigma_1$,
respectively, we get
$$
4\Sigma_0(\Sigma_0+\Sigma_1)\leq
(\Sigma_1+2\Sigma_0)^2,
$$
and that is exactly what we want.

\begin{ccc}
{\rm (Iskovskikh and Manin [IM]).} The Basic conjecture for a smooth quartic
$V\subset{\bf P}^4$
is true.
\end{ccc}

{\bf Proof.}
Obviously,
$m\leq 4n^2$.
It is a contradiction with the previous corollary.

Since it is easy to show that on
$V_4$
$|\chi|$
has no maximal cycles ([IM] or [P5]), we get

\begin{ccc}
{\rm (Iskovskikh and Manin [IM]).}  A smooth three-dimensional quartic
$V\subset{\bf P}^4$
is a birationally superrigid variety.
\end{ccc}

\section{The Sarkisov theorem on conic bundles}

We give an extremely short version of the proof of
the Sarkisov theorem [S1,2]. The idea of the proof is essentially the
same as in these well-known Sarkisov's papers. At the same time
our general viewpoint of working in codimension 1 makes the
arguments brief and very clear.

\subsection{Formulation of the theorem}

Let $S$ be a smooth projective variety,
$\mathop{\rm dim} S\geq 2,{\cal E}$
be a locally free sheaf on
$S$,
$\mathop{\rm rk} {\cal E}=3.$ Let
$$
X\subset{\bf P}({\cal E})\stackrel{\pi}{\to} S
$$
be a standard conic bundle, that is, a smooth hypersurface with
$$
\mathop{\rm Pic} X=
{\bf Z} K_X\oplus\pi^*\mathop{\rm Pic} S.
$$
Denote by
$C\subset S$
the discriminant divisor. Note that $C$ has at most
normal crossings,  the fiber over any point outside
$C$ is a smooth conic, the fiber over generic point of any component of
$C$ is a pair of distinct lines, and the inverse image of any
component of $C$ on $X$ is irreducible.

Let
$\tau:V\to F$
be another conic bundle of the same dimension
(not necessarily smooth).

\begin{ttt}
If
$|4K_S+C|\neq\emptyset$,
then any birational map
$$
\chi:X-\,-\,\to V
$$
transfers fibers into fibers, that is, there exists a map
$\bar\chi:S-\,-\,\to F$
such that
$$
\tau\circ\chi=\bar\chi\circ\pi.
$$
\end{ttt}

\subsection{Start of the proof}

Denote by
$$
{\cal F}=\{C_u|u\in U\}
$$
the proper inverse image of the family of conics
$\tau^{-1}(q), q\in F$,
and by
$$
\bar{\cal F}=\{\bar C_u=\pi(C_u)|u\in U\}
$$
its image on the
``ground''
$S$. When some birational
operations are performed with respect to these families,
the parametrizing set
$U$
is to be replaced by some dense open subset;
but for brevity we shall just keep it in mind and use the
same symbol
$U$,
meaning it to be as small as necessary.

Let
$\sigma:S^*\to S$
be a birational morphism such that:\newline
(1)
$S^*$
is projective and non-singular in codimension 1;\newline
(2) the proper inverse image
$$
{\cal F}^*=\{L_u|u\in U\}
$$
of the family
$\bar{\cal F}$
on $S^*$
is free in the following sense: for any cycle
$Z\subset S^*$
of codimension$\geq 2$
a general curve
$L_u$
does not meet
$Z$.
Existence of such a morphism
$\sigma$
can be proved quite elementary without use of the Hironaka theory
(see [P6]). Set
$$
{\bf P}^*={\bf P}(\sigma^*{\cal E}),
$$
$$
X^*=X\times_S S^*\subset{\bf P}^*,
$$
$X^*$ being a singular conic bundle over
$S^*$.
For simplicity of notations the natural morphisms of
$X^*$ to
$S^*, X$
will be denoted by
$\pi,\sigma$
respectively, and the map
$\chi\circ\sigma$
just by
$\chi$.

\begin{ppp}
There exist:
a closed subset
$Y\subset S^*$
of codimension$\geq 2$,
a nonsingular conic bundle
$$
\pi:W\to S^*\backslash Y
$$
with the non-singular discriminant divisor
$$
C^*\subset S^*\backslash Y
$$
and
$$
\mathop{\rm Pic} W\cong{\bf Z} K_W\oplus\pi^*\mathop{\rm Pic} S^*,
$$
and a fiber-wise map
$$
\lambda:X^*-\,-\,\to W,
$$
$\pi\circ\lambda=\pi.$
Moreover,
$$
|4K_{S^*}+C^*|\neq\emptyset.
$$

\end{ppp}

{\bf Proof.} We obtain
$W$
by means of fiber-wise restructuring of
$X^*$
over  the prime divisors
$T\subset S^*$
such that
$\mathop{\rm codim}\sigma(T)\geq 2$.
If $t\in{\bf C}(S^*)$
is a local equation of
$T$ on $S^*$,
then at the generic point of
$T$
the variety
$X^*$
is given by one of the two following types of equations:\newline\newline
case 1:
$x^2+t^kay^2+t^lbz^2, k\leq l,$\newline\newline
case 2:
$x^2+y^2+t^kaz^2,$\newline\newline
where
$(x:y:z)$
are homogeneous coordinates on
${\bf P}^2$,
and
$a,b$ are regular and non-vanishing at a generic point of $T$.
In the case 1 for
$k\geq 2$ the variety
$X^*$
has a whole divisor of singular points, that is,
$\pi^{-1}(T)$. Blow it up
$[k/2]$ times.
Now in both cases the singularity of our variety over $T$
is either of the type
$A_n$ or of the type $D_n$.
Blowing up the singularities,
{\it covering} $T$,
and contracting afterwards
$(-1)-$components in fibers, we get the Proposition. The last
statement is easily obtained by computing the
discrepancy of
$\nu_T$ on $S$.

Denote
$\chi\circ\lambda^{-1}$
by
$\chi:W-\,-\,\to V$
again.

Let
$Z\subset W\times V$
be the (closed) graph of
$\chi$, $\varphi$ and $\psi$
be the projections
(birational morphisms) onto
$W$ and $V$, respectively.
Obviously,
$Z$
is projective over
$W$.

\begin{ppp}
For any closed set
$Y^*\supset Y$ of codimension$\geq 2$
there exists an open set
$U\subset F$
such that
$$
\psi^{-1}\tau^{-1}(U)\subset
\varphi^{-1}\pi^{-1}(S^*\backslash Y^*)
$$
and
$\psi^{-1}\tau^{-1}(U)$
is projective over
$\tau^{-1}(U)\subset V$.
\end{ppp}

{\bf Proof:}
it follows immediately from the fact that the family of curves
${\cal F}^*$ is free on
$S^*$.

\subsection{The test surface construction}

Now let
$|H^*|$
be any linear system, which is the inverse image of a very
ample linear system on
$F$, and
$|\chi|$
be its proper inverse image on
$W$.
Write down
$$
|\chi|\subset
|-\mu K_W+\pi^*A|
$$
for some
$\mu\in{\bf Z}_+$
and
$A\in\mathop{\rm Pic} S^*$.
If
$\mu=0$,
we get the statement of the Theorem. So assume that
$\mu\geq 1$.

Let us show that this is impossible.

In the notations of the last Proposition, set
$Q=
\psi^{-1}\tau^{-1}(U)$.
Obviously, we may assume that
$$
\psi:Q\to
\tau^{-1}(U)\subset V
$$
is an isomorphism.
For a generic conic
$R_u,u\in U$,
$$
(H^*\cdot R_u)=0,
$$
$$
(K_V\cdot R_u)=-2.
$$
So the same is true on
$Q$. Thus for some prime divisors
$T_1,\dots,T_m\subset Q$
we get
$$
\left(
(-\mu\varphi^* K_W+
\varphi^*\pi^* A-
\sum^m_{i=1}a_iT_i)
\cdot
\psi^{-1}(R_u)
\right)=0,
$$
$$
\left(
(\varphi^* K_W+\sum^m_{i=1}
d_iT_i)\cdot
\psi^{-1}(R_u)
\right)=-2.
$$
Making the set
$U$
smaller if necessary, we may assume that
$$
\left(
T_i\cdot\psi^{-1}(R_u)
\right)\geq 1
$$
for all
$i$.
Thus the cycles
$$
\pi\circ\varphi(T_i)
$$
have codimension 1 in
$S^*$
and
$T_i$'s can be realized by the successions of blow-ups
$$
\begin{array}{cccc}
\displaystyle
\varphi^{(i)}_{j,j-1}: & X^{(i)}_j & \longrightarrow & X^{(i)}_{j-1} \\
\displaystyle
  & \bigcup & & \bigcup \\
\displaystyle
  & E^{(i)}_j & \longrightarrow & B^{(i)}_{j-1},
\end{array}
$$
where
$B^{(i)}_0=\varphi(T_i),
B^{(i)}_{j+1}$
covers
$B^{(i)}_j, E^{(i)}_{K(i)}=T_i.$
Since
$|\chi|$
has no fixed components,
$\mathop{\rm deg}(B^{(i)}_{j+1}\to B^{(i)}_j)=1$
and the corresponding graph of discrete valuations is a chain.
Taking the union of these blow-ups
(that is, throwing away some more cycles of codimension 2 from $S^*$),
we get on $Q$ that
$$
|\tilde\chi|\subset
\left|
-\mu\varphi^* K_W+\varphi^*\pi^* A-
\sum_{i,j}\nu_{i,j}E^{(i)}_j
\right|,
$$
whereas the canonical divisor on
$Q$
is equal to
$$
\varphi^*K_W+
\sum_{i,j}E^{(i)}_j.
$$
Consequently, as far as
$\mu\geq 1$, the divisor
$$
\varphi^*\pi^*A-
\sum_{i,j}(\nu_{i,j}-\mu)E^{(i)}_j
$$
intersects
$\psi^{-1}(R_u)$
negatively. Of course, we may assume that
$$
\nu_{i,K(i)}\geq \mu+1
$$
for all
$i=1,\dots,m$.

Now consider the surface
$\Lambda_u=\pi^{-1}(\pi\circ\varphi(\psi^{-1}(R_u))$
(the {\it test surface} -- see [P5,6])
and its proper inverse image
$\Lambda^*_u$ on $Q$. These surfaces are projective
and, since
${\cal F}^*$
is free, we get
$$
(D^2\cdot\Lambda^*)\geq 0,
$$
where
$D$
is the class of
$\psi^{-1}(|H^*|).$
On the other hand, setting
$L=\psi^{-1}(R_u),
\bar L=\pi(L)$,
we can write down
$(D^2\cdot\Lambda^*)$
as
$$
4\mu(A\cdot\bar L)-\mu^2\left(
(4K_{S^*}+C^*)\cdot
\bar L
\right)-
\sum_{i,j}\nu^2_{i,j}
(E^{(i)}_j\cdot L)
$$
(since for a generic
$u\in U$
the curve
$\psi^{-1}(R_u)$
intersects
all the $T$'s
transversally).
At the same time, according to the remark above,
$$
(A\cdot\bar L)<
\sum_{i,j}(\nu_{i,j}-\mu)(E^{(i)}_j\cdot L),
$$
so that
$$
4\mu(A\cdot\bar L)<
$$
$$
<\sum_{i,j}4\mu(\nu_{i,j}-\mu)(E^{(i)}_j\cdot L)\leq
$$
$$
\leq\sum_{i,j}\nu^2_{i,j}(E^{(i)}_j\cdot L).
$$
Since the intersection
$$
\left(
(4K_{S^*}+C^*)\cdot\bar L
\right)
$$
is obviously nonnegative, we get a contradiction:
$$
(D^2\cdot\Lambda^*)<0.
$$
Q.E.D.

\vspace{2cm}

\centerline{\Large\bf References}

\vspace{1cm}

\noindent
[N] Noether M., Ueber Fl{\" a}chen welche Schaaren rationaler
Curven besitzen.-- Math. Ann. {\bf 3}, 1871, 161-227.
\newline
\newline
[F1] Fano G., Sopra alcune varieta algebriche a tre dimensioni aventi
tutti i generi nulli.-- Atti Acc. Torino. {\bf 43}, 1908, 973-977.
\newline
\newline
[F2] Fano G., Osservazioni sopra alcune varieta non razionali aventi
tutti i generi nulli.-- Atti Acc. Torino. {\bf 50}, 1915, 1067-1072.
\newline
\newline
[F3] Fano G., Nouve ricerche sulle varieta algebriche a tre dimensioni
a curve-sezioni canoniche.-- Comm. Rent. Ac. Sci. {\bf 11}, 1947, 635-720.
\newline
\newline
[IM] Iskovskikh V.A. and Manin Yu.I., Three-dimensional quartics and
counterexamples to the L{\" u}roth problem.--
Math. USSR Sb. {\bf 15.1}, 1971, 141-166.
\newline
\newline
[I] Iskovskikh V.A., Birational automorphisms of three-dimensional
algebraic varieties.-- J. Soviet Math. {\bf 13}, 1980, 815-868.
\newline
\newline
[IP] Iskovskikh V.A. and Pukhlikov A.V., Birational automorphisms
of multi-dimensi\-onal algebraic varieties.--
J.Math.Sci. {\bf 82}, 4 (1996), 3528-3613.
\newline
\newline
[S1] Sarkisov V.G., Birational automorphisms of conical fibrations.--
Math. USSR Izv. {\bf 17}, 1981.
\newline
\newline
[S2] Sarkisov V.G., On the structure of conic bundles.--  Math. USSR
Izv. {\bf 20.2}, 1982, 354-390.
\newline
\newline
[T1] Tregub S.L., Birational automorphisms of a three-dimensional
cubic.-- Russian Math. Surveys. {\bf 39.1}, 1984, 159-160.
\newline
\newline
[T2] Tregub S.L., Construction of a birational isomorphism of a
three-dimensional cubic and a Fano variety of the first kind with
$g=8$, connected with a rational normal curve of degree 4.--
Moscow Univ. Math. Bull. {\bf 40.6}, 1985, 78-80.
\newline
\newline
[Kh] Khashin S.I., Birational automorphisms of a double Veronese
cone of dimension 3. -- Moscow Univ. Math. Bull. {\bf ?.1}, 1984, 13-16.
\newline
\newline
[P1] Pukhlikov A.V., Birational isomorphisms of four-dimensional
quintics.-- Invent. Math. {\bf 87}, 1987, 303-329.
\newline
\newline
[P2] Pukhlikov A.V., Birational automorphisms of a double space and a
double quadric.-- Math. USSR Izv. {\bf 32}, 1989, 233-243.
\newline
\newline
[P3] Pukhlikov A.V., Birational automorphisms of a three-dimensional
quartic with an elementary singularity.-- Math. USSR Sb. {\bf 63}, 1989,
457-482.
\newline
\newline
[P4] Pukhlikov A.V., Maximal singularities on a Fano variety $V^3_6$.--
Moscow Univ. Math. Bull. {\bf 44}, 1989, 70-75.
\newline
\newline
[P5] Pukhlikov A.V., A note on the theorem of V.A.Iskovskikh and
Yu.I.Manin on the three-dimensional quartic.-- Proc. Steklov Math. Inst.
{\bf 208}, 1995.
\newline
\newline
[P6] Pukhlikov A.V., Birational automorphisms of singular double
spaces.-- Cont. Math. and Appl. 1995, {\bf 11}.
\newline
\newline
[R] Reid M., Birational geometry of 3-folds according to
Sarkisov.
\newline
\newline
[C] Corti A., Factoring birational maps of threefolds
after Sarkisov.
J. Algebraic Geometry {\bf 4} (1995) 223--254.
\newline
\newline
[K] Koll{\' a}r J., Nonrational hypersurfaces.
J. of Amer. Math. Soc. {\bf 8} (1995) 241--249.

\end{document}